\documentclass[11pt]{amsart}
\usepackage{amsmath}
\usepackage{xy}
\usepackage{amscd}
\newfont{\sheaf}{eusm10 scaled\magstep1}
\def\hol{{\mathcal{O}}}
\newtheorem{definition}{Definition}[section]

\newtheorem{proposition}{Proposition}[section]
\newtheorem{corollary}[proposition]{Corollary}
\newtheorem{lemma}[proposition]{Lemma}
\newtheorem{theorem}[proposition]{Theorem}

\newtheorem{remark}{Remark}[section]

\DeclareMathOperator{\Hilb}{Hilb}
\DeclareMathOperator{\PGL}{PGL}
\DeclareMathOperator{\Ker}{Ker}

\DeclareMathOperator{\Sing}{Sing}
\DeclareMathOperator{\Aut}{Aut}

\DeclareMathOperator{\Pic}{Pic}

\title[Rationality of $\mathcal{R}_{4,\langle 3 \rangle}$] {The rationality of the moduli space of genus four curves endowed with an  
order three subgroup of their Jacobian}
\thanks{Part of the article was written while the second author was visiting professor in Bayreuth financed by the DFG-Forschergruppe 790 {\it Classification of algebraic Surfaces and compact complex Manifolds} and by the research program PRIN 2006-08: {\it Geometria delle variet\'a algebriche e dei loro spazi di moduli}}
\author{Ingrid Bauer}
\address{Mathematisches Institut, Universit\"at Bayreuth, NW II, D-95440 Bayreuth; email:
Ingrid.Bauer@uni-bayreuth.de}
\author{Alessandro Verra}
\address{Universita' di Roma III, Dipartimento di Matematica, L.go S. Leonardo Murialdo, 1 00146 Roma;  email: verra@mat.uniroma3.it}

\begin{document}
\maketitle

\section*{Introduction} \noindent
Let $C$ be 
a smooth, irreducible complex projective curve of genus $g$ and let $\eta \in \Pic^{0}(C)$ be a non trivial $n$-th root of the trivial bundle $\mathcal O_{C}$. For several different reasons a special
attention has been payed, now and in the recent past, to the moduli spaces 
$
\mathcal R_{g,n}
$
of pairs $(C,\eta)$ as above, and to its possible compactifications, (see e.g., \cite{ccc}). \par \noindent
For instance, they are generalizations of the case $n = 2$, the so called \it Prym moduli spaces, \rm usually denoted
by $\mathcal R_{g}$. Since they are related to the theory of Prym varieties, the interest in this case occupies a prominent
position. In particular, many results on the Kodaira dimension of $\mathcal R_{g}$ are now available, while classical geometric descriptions of $\mathcal R_{g}$ exist for $g \leq 7$. More precisely, let us mention that G. Farkas and K. Ludwig proved that $\mathcal R_{g}$ is of general type for $g \geq 14$ and $g \neq 15$, (\cite{farkasludwig}).  On the other hand, unirational parametrizations of $\mathcal R_{g}$ are known for $g \leq 7$, (\cite{fabrat}, \cite{dolga}, \cite{ils}, \cite{donagi}, \cite{verra84}, \cite{verra08}). \par \noindent
One can also consider the moduli spaces $\mathcal R_{g,\langle n \rangle}$ of pairs $(C, \mathbb Z/n \mathbb Z)$, where $C$ is a smooth, irreducible complex projective curve of genus $g$ and $\mathbb Z/n \mathbb Z$ is a cyclic subgroup of order $n$ of $\Pic^{0}(C)$. As $\mathcal R_{g,\langle 2 \rangle} = \mathcal R_{g, 2}$, these mouli spaces are generalizing the Prym moduli spaces in a (slightly) different way.
In contrast to the case $n=2$, not very much is known about $\mathcal R_{g,n}$ resp. $\mathcal R_{g,\langle n \rangle}$ for $n > 2$. In particular, the (probably short)  list of all pairs $(g,n)$
such that $\mathcal R_{g,n}$ resp. $\mathcal R_{g,\langle n \rangle}$ has negative Kodaira dimension is not known. \par \noindent The rationality of $\mathcal R_{g,n}$ resp. $\mathcal R_{g,\langle n \rangle}$ has been proved in some cases of very low genus:  the case of $\mathcal R_{4}$ is a  result of Catanese, (\cite{fabrat}). The rationality of $\mathcal R_{3}$ was proved by Katsylo in \cite{katsylo}. Independent proofs are also due to Catanese and to Dolgachev: see \cite{dolga} (also for $\mathcal R_{2}$). Recently,  the rationality of $\mathcal R_{3,3}$ and of $\mathcal R_{3,\langle 3 \rangle}$ has been proven by Catanese and the first author (cfr. \cite{bauercat}).
\par \noindent 
To complete the picture, we recall that $\mathcal R_{1,n}$ is an irreducible curve for every prime $n$
and that its geometric genus is well known. \par  \noindent Usually, for reaching one of the previous rationality results, beautiful geometric properties of low genus canonical curves come to the rescue. 

\medskip\noindent
This is the case also for  the present note:  we will  use the classical geometry of cubic surfaces in $\mathbb P^{3}$ to prove the following
\begin {theorem} The moduli space $\mathcal R_{4,\langle 3 \rangle}$ is rational.
\end{theorem} \rm \par \noindent
Let us very briefly  describe our  proof.  \medskip \par  \noindent
We consider the moduli space $\mathcal P$ of the sets $\underline x $ of six  general (unordered) points in $\mathbb P^{2}$. Equivalently,  $\mathcal P$ is the moduli space of pairs $$(S, \sigma),$$ such that  $S \subset \mathbb P^{3}$ is a smooth cubic surface and
$\sigma: S \to \mathbb P^{2}$
is the blow up of $\underline x$.  Let $C \in | \omega_{S}^{-2} |$ be smooth,  and let $L \in | \sigma^{*} \mathcal O_{\mathbb P^{2}}(1)|$. It turns out that $C$ is a canonical curve of genus 4 endowed with the line bundle
$$
\eta_{C} := \omega_{C}(-L) \in \Pic^{0}(C).
$$
We say that $(C, \eta_{C})$ defines a point of $\mathcal R_{4,3}$ if $C$ is smooth, $\eta_{C}$ is non trivial and $\eta_{C}^{3} \cong \mathcal O_{C}$. Of course, we are interested in the locally closed set 
$$
Q_{\underline x} := \lbrace C \in  | \omega_{S}^{-2} | \ | \  (C, \eta_{C}) \text { defines  a point of}  \ \mathcal R_{4,3} \rbrace.
$$
To study this set, (as in  \cite{barth}) we make use of the cup product map
$$
\mu: H^{1}({\omega_{S}}^{-1}(-3L)) \otimes H^{0}({\omega_{S}}^{-2}) \to H^{1}(\omega_{S}^{-3}(-3L)).
$$
Let $v \otimes s \in H^{1}({\omega_{S}}^{-1}(-3L)) \otimes H^{0}({\omega_{S}}^{-2})$ be a non zero decomposable vector and $C = div(s)$: we show that $v \otimes s \in \Ker \mu$ iff $(C, \eta_{C})$ defines a point of $\mathcal R_{4,3}$. Then we  consider the Segre product
$$
\Sigma :=  \mathbb P(H^{1}(\omega_{S}^{-1}(-3L))) \times | \omega_{S}^{-2}| \ \subset \ \mathbb P := \mathbb P (H^{1}({\omega_{S}}^{-1}(-3L)) \otimes H^{0}({\omega_{S}}^{-2}))$$ and the  intersection scheme
$$
\mathbb M_{\underline x}^{o} := \mathbb P(\Ker  \mu) \cap \Sigma^{o} \subset \mathbb P, $$ where 
$\Sigma^{o} := \lbrace (z,C) \ \in \ \Sigma \ | \ C \ \rm is \ smooth \rbrace. $ It turns out that the projection $(z,C) \to C$ induces a biregular map
$$
q_{\underline x}: \mathbb M_{\underline x}^{o} \to Q_{\underline x}.
$$
Moreover, by the previous remarks, there exists a natural rational map
$$
\phi_{\underline x}: \mathbb M_{\underline x}^{o} \to \mathcal R_{4,3},
$$
sending $(z,C)$ to the moduli point of the pair $(C, \eta_{C})$. \bigskip  \par \noindent
As a first step, we show that $\mathbb M_{\underline x}^{o}$  is integral of the expected dimension $5$. Let $\mathbb M_{\underline x}$ be the Zariski closure of $\mathbb M_{\underline x}^{o}$. Then we  show that the projection 
$$
p_{\underline x}: \mathbb M_{\underline x} \to \mathbb P(H^{1}(\omega_{S}^{-1}(-3L)))
$$
is a locally trivial $\mathbb P^{1}$-bundle, over a suitable non empty open set. \bigskip  \par \noindent
As a second step, we globalize the construction to the moduli space $\mathcal P$. On a  dense open set  $\mathcal P^{o} \subset \mathcal P$, we define vector bundles $\mathcal E$ and $\mathcal H$ such that: \\
- the fibre of $\mathcal E$ at the moduli point of $\underline x$ is $H^{1}(\omega_{S}^{-1}(-3L))$, \\
- the fibre of $\mathcal H$ at the moduli point of $\underline x$ is $H^{0}(\omega_{S}^{-2})$. \\
Let
$
\mathbb E := \mathbb P (\mathcal E) \ \rm and \ \mathbb H := \mathbb P (\mathcal H)
$
be their projectivizations. In the fibre product $\mathbb E \times_{\mathcal P^{o}} \mathbb H$ we  construct then an  integral variety
$$
\mathbb M \subset \mathbb E \times_{\mathcal P^{o}} \mathbb H
$$
such that, at the moduli point of $\underline x$, the fibre of the projection $\pi: \mathbb M \to \mathcal P^{o}$ is  the Zariski closure $\mathbb M_{\underline x}$  of $\mathbb M_{\underline x}^{o}$. 
This allows us to define a rational map
$$
\phi: \mathbb M \to \mathcal R_{4,3},
$$
whose restriction to the fibre $\mathbb M_{\underline x}$ is the above considered rational map $\phi_{\underline x}$. \bigskip \par \noindent
As a third step we show that: \medskip \par \noindent \it
(1) $\phi: \mathbb M \to \mathcal R_{4,3}$ is birational. \medskip  \par \noindent
(2) $\mathbb M$ is birational to $\mathcal P \times \mathbb P^{4} \times \mathbb P^{1}$.
\medskip \par \noindent \rm
Indeed, (2) follows because $\mathbb E$ is a $\mathbb P^{4}$-bundle over $\mathcal P^{o}$ and because  $\mathbb M$ is birational to a $\mathbb P^{1}$-bundle over $\mathbb E$. \\
Finally, we observe that $\mathcal R_{4,\langle 3 \rangle} = \mathcal R_{4,3} /< i>$, where $i: \mathcal R_{4,3} \rightarrow \mathcal R_{4,3}$ is the involution mapping $(C, \eta)$ to $(C, \eta^{-1})$. Moreover this involution is given on $\mathcal P \times \mathbb P^{5}$ by the Schlaefli involution $j: \mathcal P \to \mathcal P$ on the first factor and the identity on the second factor. The quotient $\mathcal P/ <j>$ turns out to be the moduli space of a double six, a well known configuration of lines in a cubic surface, (see section 2). Now, the rationality of $\mathcal R_{4,\langle 3 \rangle}$ is a consequence of the following
\begin{theorem}[I. Dolgachev] The moduli space $\mathcal P /< j>$ of double sixes in $\mathbb P^{2}$ is rational. \end{theorem} \rm \par \noindent
The proof of this  theorem is an unpublished result of Igor Dolgachev.  We reproduce here his proof:  we would like to thank him warmly for allowing us to include his proof and for some useful remarks  
on this paper. \par \noindent This work also profited of some interesting  discussions on the subject with Fabrizio Catanese, during a visit of the second author in  Bayreuth. Finally we wish to thank the referee for some useful remarks.

\section{Plane sextics of genus 4} \noindent
Let $C$ be a smooth, irreducible, complete curve of genus $g \geq 3$ and let $\eta \in \Pic^{0}(C)$ be 
a non trivial line bundle. We are interested in  the linear system $|\omega_{C} \otimes\eta^{-1}|$ and its associated rational map
$$
\phi_{\eta}: C \to \mathbb P^{g-2}.
$$ 
We first derive some basic properties of $|\omega_{C} \otimes\eta^{-1}|$.
\begin{proposition}\label{basicprop} The following conditions are equivalent: 
\begin{itemize}
\item[(1)] $|\omega_{C} \otimes\eta^{-1}|$ has a base point $p$; 
\item[(2)] $h^{0}(\eta(p)) = 1$; 
\item[(3)] $\eta \cong \mathcal O_{C}(q-p)$ for some $q \in C \setminus \{p\}$.
\end{itemize}
\end{proposition}
\begin{proof} (1) $\Leftrightarrow$ (2): \ $p$ is a base point of  $| \omega_{C} \otimes\eta^{-1}|$ $\Leftrightarrow$ $ \dim |(\omega_{C} \otimes\eta^{-1})(-p)| =$ $\dim \ |\omega_{C} \otimes\eta^{-1}| $ $=$ $g-2 \Leftrightarrow h^{1}(\eta(p)) = g-1$ $\Leftrightarrow$ $h^{0}(\eta(p)) = 1$. \\
(2) $\Leftrightarrow$ (3): obvious. \end{proof} \rm \par \noindent
The next (well known) statement complements the previous one:

\begin{proposition} Let $b$ be the base divisor of $| \omega_{C}(p - q)|$, where $p,q \in C$ and $p \neq q$.
Then either $b = p$ or $C$ is hyperelliptic and $b = p + i(q)$, $i$ being the hyperelliptic involution.
\end{proposition} 

\noindent
From now on we will restrict ourselves to the case: 
$$
g = 4, \ \ \eta^{\otimes 3} \cong \mathcal O_{C}.
$$

\begin {proposition} Let $C$ be a general curve of genus $4$ and $\eta$ a non zero $3$- torsion element of
$\Pic^{0}(C)$. Then $|\omega_{C} \otimes\eta^{-1}|$ has no base points. 
\end{proposition}
\begin{proof} By prop. (\ref{basicprop}) $|\omega_{C} \otimes\eta^{-1}|$ has a base point $\Leftrightarrow$ $\eta \cong \mathcal O_{C}(q-p)$, with $p, q \in C$. In particular, $3p \sim 3q$. Let $f: C \to \mathbb P^{1}$ be the map defined by  $|3p| = |3q|$: then  $f$  has double ramification in $p$ and $q$. By Hurwitz's formula, the degree of the ramification divisor $R$ of $f$ is $12$. Then, since $2p+2q \leq R$, $f$ has at most $10$ branch points. Hence, by Riemann's existence theorem, the map $f$ depends on 7 moduli at most.  This implies the statement. \end{proof}
 
\noindent
Therefore we will assume from now on that {\em $| \omega_{C} \otimes \eta^{-1} |$ has no base points}. 

\noindent
Consider
$$
\Gamma_{\eta} \ := \ \phi_{\eta}(C) \ \subset \ \mathbb P^{2}.
$$
Since $\phi_{\eta}$ is a morphism, we have the following possibilities: 
\begin{itemize}
\item $\phi_{\eta}: C \to \Gamma_{\eta}$ has degree 1 and $\Gamma_{\eta}$ is an integral sextic, 
\item $\phi_{\eta}: C \to \Gamma_{\eta}$ has degree 2 and $\Gamma_{\eta}$ is an integral cubic, 
\item $\phi_{\eta}: C \to \Gamma_{\eta}$ has degree 3 and $\Gamma_{\eta}$ is an integral conic.
\end{itemize} \noindent
All these cases actually occur, but we will show that for a general curve $C$,  $\Gamma_{\eta}$ is a sextic and  $\Sing \Gamma_{\eta}$ consists of $6$ ordinary double points  in general position. 

\medskip
\noindent
To this purpose we consider the image $S \subset \Pic^{2}(C)$ of the Abel map
$$
a: C^{(2)} \to \Pic^{2}(C).
$$
Then the cohomology class of $S$ is $\frac {\theta^{2}}2$, where $\theta$ is the class of a theta divisor in $\Pic^{2}(C)$. In particular, we have $S^{2}= 6$. For any $\eta \in \Pic^{0}(C)$, let  $S_{\eta}$ be the translate of $S$ by $\eta$. Moreover, let 
$$
Z_{\eta} \ := \ a^{*}S_{\eta}
$$
be the pull back of $S_{\eta}$ under the Abel map. By the transversality of a general translate, $S_{\eta}$ is transversal to $a$ for general $\eta$. Therefore, in this case, the scheme $Z_{\eta}$ consists of six, smooth and distinct, points. 

\begin{lemma}\label{dinZeta} For $d \in C^{(2)}$ the following conditions are equivalent: 
\begin{itemize}
 \item[(1)] $\phi_{\eta}$ contracts $d$ to a point, 
 \item[(2)] $h^{0}(\eta(d)) \geq 1$, 
\item[(3)] $d \in Z_{\eta}$.
\end{itemize}
 \end{lemma}
\begin{proof} By assumption, $|\omega_{C} \otimes\eta^{-1}|$ has no base points. Then $\phi_{\eta}$ contracts $d$  to a point iff  $h^{0}((\omega_{C} \otimes\eta^{-1})(-d)) \geq 2$. On the other hand, $h^{0}((\omega_{C} \otimes\eta^{-1})(-d)) \geq 2$ $\Leftrightarrow$ $h^{0}(\eta(d)) \geq 1$ $\Leftrightarrow$ $d \in Z_{\eta}$. \end{proof} 
\noindent
In the following lemma we prove that, with respect to $Z_{\eta}$ and $\phi_{\eta}$, a general hyperelliptic curve $C$ 
has a sufficiently general behaviour:

\begin{lemma}\label{birat} Let $C$ be a general hyperelliptic curve of genus $4$ and let $\eta$ be
any non trivial line bundle on $C$ such that $\eta^{3} \cong \mathcal O_{C}$. Then:
\begin{itemize}
 \item[(1)] each element  $d \in Z_{\eta}$ is a smooth divisor of degree $2$; 
 \item[(2)] $\phi_{\eta}: C \to \Gamma_{\eta}$ is a birational morphism;
 \item[(3)] $Z_{\eta}$ is a smooth, 0-dimensional scheme supported in $6$ points;
 \item[(4)] $\Gamma_{\eta}$ is a sextic with an ordinary singular point of multiplicity $4$.
\end{itemize}
Moreover, conditions (1), (2), (3) hold also for a general $C$ of genus $4$.   
\end{lemma}
\begin{proof} (1) Assume by contradiction that $d \in Z_{\eta}$ is not smooth. Then it follows that $d = 2p$,
for some $p \in C$. On the other hand, lemma (\ref{dinZeta}) implies that  $\eta(d) \cong \mathcal O_{C} (d')$, where $d'$ is an effective divisor. Therefore we have $d' = q + r$, for some $q,r \in C$. Since $\eta^{3}$ is trivial the
divisors $6p$ and $3q + 3r$ generate a pencil $P$. To prove (1), we show that a pencil like $P$ cannot exist on a general hyperelliptic curve $C$. \par \noindent
First, let us show that $P$ has no base points: if   $P$ has a base point $x$,  then $x=p$ and $p \in \lbrace q, r \rbrace$. Hence $3p$ is a fixed component of $P$. Moreover $|3p|$ is a pencil of degree $3$ and, since $C$ is hyperelliptic, it has a base  point. But then, $4p$ is a fixed component of $P$ and $p = q = r$: a contradiction. \par \noindent
Now let $|h|$ be the hyperelliptic pencil of $C$. We consider the morphism
$
\psi: C \to \mathbb P^{1} \times \mathbb P^{1},
$
defined by $P \times |h|$, and the curve $\Gamma := \psi(C)$. Two cases are possible:
\begin{itemize}
\item[(A)] $\Gamma$ has bidegree $(6,2)$ and $\psi: C \to \Gamma$ is birational, 
\item[(B)] $\Gamma$ has bidegree $(3,1 )$ and $\psi: C \to \Gamma$ has degree two. 
 \end{itemize}
In both cases we show that $C$ cannot be a general hyperelliptic curve. \medskip \par \noindent
(A) Note that $\psi(q)$ and $\psi(r)$ are in the same line $l$ of type $(1,0)$ and $\psi(p) \notin l$. Indeed the pull-back by $\psi$ of the ruling of lines of type $(1,0)$ is the pencil $P$, generated by the divisors $6p$ and $3q + 3r$. \\ Now let us fix coordinates $(X_{0}:X_{1}) \times (Y_{0}:Y_{1})$ on $\mathbb P^1 \times \mathbb P^1$, such that
$$ l =  \{ Y_{1} = 0 \} \ , \ \psi(p) = \{ X_0 = Y_0 = 0 \} \ , \ \lbrace \psi(q) , \psi(r) \rbrace = \{ X_0X_e = Y_1 =  0 \},$$
where $e = 1$, if $\psi(q) \neq \psi(r)$, and $e = 0$, if $\psi(q) = \psi(r)$. After this choice $\Gamma$
is an element of the $8$-dimensional  linear system $\Sigma$ defined by the equation
$$
aY_{0}^{2}X_0^3X_e^3 + bY_{1}^{2}X_{0}^{6} +  \sum_{i=0 \dots 6}c_{i}Y_{0}Y_{1}X_{0}^{i}X_{1}^{6-i}= 0.
$$
A general element of $\Sigma$ is smooth, moreover $\dim \Sigma = 8$. On the other hand, $p_{a}(\Gamma) = 5$, whence $\Gamma$ is singular and belongs to the $7$-dimensional discriminant hypersurface $\Delta \subset \Sigma$. Note that the stabilizer $G \subset \Aut (\mathbb P^1 \times \mathbb P^1)$ of $\Sigma$ has dimension $\geq 1$, whence $\dim \Delta / G \leq 6$. But then, $C$ has at most $6$ moduli and it is not general hyperelliptic.
\medskip \par \noindent
(B) $\Gamma$ is a smooth rational curve of type $(3,1)$. Note that  $6p = \psi^*l_{1}$  and $3q + 3r = \psi^{*}l_{2}$, where $l_{1}, l_{2}$ are lines of type $(0,1)$. The latter equality
implies $l_{2} \cdot \Gamma = 3v$ for some $v \in l_{2}$, hence $q+r = \psi^{*}v$. The first one implies $l_{1} \cdot \Gamma = 3u$, for some $u \in l_{1}$, hence $2p = \psi^{*}u$. Since $\Gamma$ is rational, it follows that $\psi^{*}(u-v) = 2p - (q+r) \sim 0$. Then $\eta \cong \mathcal O_{C}$ and this case
is excluded. 
\medskip \par \noindent
(2) A non birational $\phi_{\eta}$ ramifies at some $p \in C$. Hence $\phi_{\eta}$ contracts $2p$ and, by lemma (\ref{dinZeta}), $2p \in Z_{\eta}$. By (1) this is impossible for a general $C$.  \medskip \par \noindent
(4) Let $|h|$ be the $g^1_2$ of $C$, then $\omega_{C} \cong \mathcal O_{C}(3h)$. This implies that
$|\omega_{C} \otimes \eta^{-1}|  = |h + b|$, where $b \in  | \eta^{-1}(2h) |$. $b$ has degree 4 and $\phi_{\eta}(b)$ is a point $o$. Then, since $\phi_{\eta}$  is birational, $o$ is a $4$-uple point of $\Gamma_{\eta}$. By the genus formula, $\Sing \Gamma_{\eta} = \{ o \}$ and, by (1), $b$ is smooth. Hence  $o$ is an ordinary $4$-uple point. 
\medskip \par \noindent
(3) Let $b = p_1 + p_2 + p_3 + p_4$ be the divisor considered above, then the conclusion is definitely clear:   by  (4)  $Z_{\eta}$ is a scheme supported in the $6$ points \ 
$
p_{i} + p_{j}, \ \ 1 \leq i < j \leq 4.
$ \
On the other hand, we know that $Z_{\eta}$ has length 6. So  $Z_{\eta}$ is smooth. 

\medskip\noindent
Finally we remark that  statements (1), (2), (3) define an open subset $U$, in the moduli space of curves of genus 4, where the statements are true. By the previous part of the  proof, $U$ is not empty. This completes the proof. 
\end{proof}
\begin{lemma}\label{nodes} If $C$ is general, $\Gamma_{\eta}$ has no point of multiplicity $ m \geq 3$.
\end{lemma} 
\begin{proof}  We can assume that $\Gamma_{\eta}$ is a non hyperelliptic sextic. Then each point
$o \in \Gamma_{\eta}$ has multiplicity $m \leq 3$. To exclude the case $m=3$, we consider the  theta divisor
$\Theta  := \{N \in \Pic^{3}(C) \ | \ h^{0}(N) \geq 1 \}$ and observe: \noindent
{\bf Claim.} $\Gamma_{\eta}$ has a triple point $\Leftrightarrow$   $L \otimes \eta \in \Theta$  $\Leftrightarrow $ $M \otimes \eta^{-1} \in \Theta$, where $|L|$ and $|M|$ are the only trigonal pencils on $C$. \par \noindent
{\em Proof of the claim.}
Note that the following conditions are equivalent: \\ (1) $\Gamma_{\eta}$ has a triple point $o$, \\ (2)   there exists  an effective $d \in Div^3(C)$ such that $\phi_{\eta}(d) = o$, \\ (3)  there exists an effective $d \in Div^3(C)$ such that $ h^0(\omega_C \otimes \eta^{-1}(-d)) = 2$, \\ (4) there exists an effective  $d \in Div^3(C)$  such that $h^0( \eta(d)) = 2$. \\
On the other hand a non hyperelliptic $C$ has at most two trigonal line bundle, say $L$ and $M$. Therefore (2) holds iff $\lbrace L, M \rbrace$ $=$ $\lbrace \omega_C \otimes \eta^{-1}(-d), \eta(d) \rbrace$ iff $L \otimes \eta \in \Theta$ iff  $M \otimes \eta^{-1} \in \Theta$. This proves the claim. \\
The previous conditions cannot hold for all $\eta \in Pic^0_3(C) - \lbrace \mathcal O_C \rbrace$. Indeed  this would imply that  the set $L + Pic^0_3(C) := \lbrace \eta \otimes L, \ \eta \in Pic^0_3(C) \rbrace$ is in the theta divisor $\Theta$ of $Pic^3(C)$. Equivalently $Pic^0_3(C)$ would be contained in  $$ \Theta_0 :=  \lbrace N \otimes L^{-1}, \  N  \in \Theta \rbrace = \lbrace  N^{-1} \otimes M, \ N \in \Theta \rbrace, $$
where the latter equality follows from $N \in \Theta \Leftrightarrow \omega_C \otimes N^{-1} \in \Theta$ and $\omega_C \cong L \otimes M$. But it is well known that, in the embedding defined by $3\Theta_0$, the set $Pic^0_3(C) - \lbrace \mathcal O_C \rbrace$  is not  in a hyperplane. Hence there 
exists a non trivial $\eta$ such that $h^0(\eta(d)) \leq 1$ for each effective $d \in Div^3(C)$. Since $\mathcal R_{4,3}$ is irreducible, this property holds $\forall \  \eta \in Pic^0_3(C) - \lbrace \mathcal O_C \rbrace$ if $C$ is general.  \end{proof}
 \noindent
The next definition is well known: it will be  important in the sequel.
\begin{definition} Six distinct points of $ \mathbb P^{2}$ are in general position if no conic contains all of them and no line contains three of them.  \end{definition} \rm \par \noindent
We are now ready to prove the main result of this section:
\begin{theorem}\label{ingenpos}  Let $C$ be a general curve of genus $4$ and let $\eta$ be
any non trivial line bundle on $C$ such that $\eta^{3} \cong \mathcal O_{C}$. Then: \\
(1) $\phi_{\eta}: C \to \Gamma_{\eta}$ is a birational morphism, \\
(2) $\Sing \Gamma_{\eta}$ consists of six ordinary nodes in general position.
\end{theorem}
\begin{proof} From the previous lemmas we know that $	\phi_{\eta}$ is a birational morphism and that $\Gamma_{\eta}$ is a  sextic with at most double points. Let $Z_{\eta}$ be the scheme considered in
(\ref{dinZeta}). We know from lemma (\ref{birat}) that $Z_{\eta}$ consists of six, distinct and smooth, effective divisors of degree $2$ on $C$, and that each of them is contracted by $\phi_{\eta}$ to a double point of $\Gamma_{\eta}$. It is  easy to deduce that $\Sing \Gamma_{\eta}$ consists of six ordinary nodes: $x_{1}, \dots, x_{6}$.   It remains to show that they are in general position. \\
Assume that a conic $B$ contains $\Sing \Gamma_{\eta}$. It is easy to deduce that, for any line $L \subset \mathbb P^{2}$, $\phi^{*}L \in |\omega_C|$, whence $|\omega_C| = |\omega_C \otimes \eta^{-1}|$: a contradiction. \\
The condition that no line contains $3$ points of $\Sing  \Gamma_{\eta}$ is equivalent to the following condition on $Z_{\eta}$: for any three distinct elements $u, v, t \in Z_{\eta}$ their sum $u+v+t$ is not
in $|\omega_{C} \otimes \eta^{-1}|$. So it suffices to prove the latter condition for at least one pair $(C,\eta)$. 
Let $C$ be general hyperelliptic: from the proof of lemma (\ref{birat}), (4) we know that 
$$
Z_{\eta} = \lbrace x_{i} + x_{j} \ 1 \leq i < j \leq 4 \rbrace, \ \rm for \  some \ b = x_{1}+x_{2}+x_{3}+x_{4}, 
$$
and  that $|\omega_{C} \otimes \eta^{-1}| = |b + h|$, where $|h|$ is the hyperelliptic pencil of $C$. Let $u+v+t \in |\omega_{C} \otimes  \eta^{-1}|$ for some $u,v,t \in Z_{\eta}$, then $u + v + t \sim b + h$. Moreover $u + v + t - b$ is effective because $ Supp \ b \subset Supp \ u \cup Supp \ v \cup Supp \ t$. Hence $u + v + t - b = p' + p''$, where $p',p'' \in Supp \ b$ and  $b \equiv h + p' + p''$. But then, $|\omega_C \otimes \eta^{-1}| = |2h+p'+p''|$ and $p', p''$ are base points: a contradiction.  \end{proof} \rm \par \noindent
We end this section giving some definitions and fixing some further notations, which will be used in the subsequent chapters. \par \noindent
Let $\Hilb_{6}(\mathbb P^{2})$
be the Hilbert scheme of 6 points in $\mathbb P^{2}$, then:
\begin{definition} $\mathcal X$ is the open subset of $\Hilb_{6}(\mathbb P^{2})$ parametrizing those
schemes $\underline x$ which are supported in six points in general position.
\end{definition}
\begin{definition} 1) $\mathcal R_{4,3}$ is the moduli space of pairs $(C,\eta)$ such that: 
\begin{itemize}
 \item[(1)] $C$ is a smooth, irreducible projective curve of genus 4,
 \item[(2)] $\eta$ is a non zero 3-torsion element of $\Pic^{0}(C)$.
\end{itemize}
\end{definition} \rm \par \noindent
$\mathcal R_{4,3}$ is an irreducible, quasi projective variety of dimension $9$. Throughout the paper we will keep the previous notations.

\begin{remark} \rm
Note that we have a natural involution $i$ on $\mathcal R_{4,3}$, sending $(C,\eta)$ to $(C,\eta^{-1})$. We will denote the quotient $\mathcal R_{4,3} / i$ by this involution by $\mathcal R_{4,\langle 3 \rangle}$. Obviously, $\mathcal R_{4,\langle 3 \rangle}$ is the moduli space of pairs of smooth, irreducible projective curves of genus $4$ together with a cyclic subgroup of $\Pic^0(C)$ of order $3$.
\end{remark}
\section {A cup product  map on the cubic surface} 

\noindent
It is important to point out that the previous theorem (\ref{ingenpos}) relates the family of pairs  $(C,\eta)$ as above
to \it smooth \rm cubic surfaces in $\mathbb P^{3}$. \\ More precisely, assume that $(C, \eta)$ defines  a general point of $\mathcal R_{4,3}$. Then 
$
\Sing \Gamma_{\eta}
$
is an element of $\mathcal X$, i.e., its points are in \it general position. \rm Let 
$
\sigma: S \to \mathbb P^{2}
$
be the  blow up of  $\Sing  \Gamma_{\eta}$. Then $S$ is a Del Pezzo surface of degree 3, i.e.,  its anticanonical divisor $-K_{S}$ is \it semiample \rm and $K_{S}^{2} = 3$. It is  well known that, for the  blow up of $\mathbb P^{2}$ in six distinct points,  the following conditions are equivalent:  \it
\begin{itemize}
 \item the six points are in general position,
 \item the anticanonical divisor is very ample.
\end{itemize} \rm \par \noindent
  Thus we will assume that our  \it $S$ is anticanonically embedded in $\mathbb P^{3}$ as a smooth cubic surface. \rm  We also remark that:
\begin{itemize}
 \item[(1)] $C \in |-2K_{S}|$,
 \item[(2)] $\eta \cong \mathcal O_{C}(-K_{S}-L)$,
 \item[(3)] $E \cdot C = \sum_{d \in Z_{\eta}} d$,
\end{itemize} 

\noindent
where $L \cong \sigma^{*}\mathcal O_{\mathbb P^{2}}(1)$ and $E$ is the exceptional divisor  of $\sigma$. 
In particular, $C$ is a quadratic section of $S$ and a canonical curve of genus four. 

\begin{remark} \rm
Though $\Sing \Gamma_{\eta}$ is a set of points in general position, still we did not prove that it is
a {\em general point} of $\mathcal X$, so that $S$ is a {\em general} smooth cubic surface. The proof of this property is  a relevant step of this section (see theorem (\ref{fibre})). 
\end{remark}
In the sequel we  will deal with \it any element \rm
$$
\underline x \in \mathcal X,
$$
supported on the set $\lbrace x_{1} \dots x_{6} \rbrace$. Let $ \sigma: S \to \mathbb P^{2}$
be the  blow up of $\underline x$ and let $E_{i} = \sigma^{-1}(x_{i})$, $i = 1 \dots 6$. As usual, we will have the line bundle
$$
L := \sigma^{*}\mathcal O_{\mathbb P^{2}}(1)
$$
and the \it exceptional divisor \rm $ E := E_1 + \dots + E_6$ of $\sigma$.
\begin{definition} For any $C \in | -2K_{S} |$ we define:
 \begin{itemize}
\item $s_C :=$ any non zero vector of $H^0(\mathcal O_{S}(-2K_S))$ vanishing on $C$,
\item  $\eta_C := \mathcal O_C(-K_S-L) $,
\item $n_{C} := C \cdot E$.
\end{itemize}
\end{definition} 
\noindent
It is easily checked that $\eta^{\otimes 3}(n_{C}) \cong \mathcal O_{C}(-2K_{S})$.
\begin{proposition}  For any  $C \in |-2K_{S}|$ the sheaf  $\eta_{C}$ is non trivial. 
\end{proposition}
\begin{proof} Consider the long exact sequence associated to the exact sequence
$$
0 \to \mathcal O_S(K_S - L) \to \mathcal O_S(-K_S-L) \to \eta_C \to 0.
$$
We have $h^1(\hol_S(K_S - L)) = 0$. Since no conic contains $\underline x$, we have also $h^{0}(\mathcal O_{S}(-K_{S}-L)) = 0$. Then $h^{0}(\eta_{C}) = 0$ and $\eta_{C}$ is non trivial. \end{proof}  
\noindent
Note that $\sigma |C$ is the map defined by $| \omega_C \otimes \eta^{-1}_C |$. Therefore, if $(C,\eta)$ is a pair as in section 1 and $\underline x = \Sing \Gamma_{\eta}$, one has $\eta \cong \eta_{C}$.
\bigskip \par \noindent
We  want to understand  the family of smooth curves $C \in |-2K_{S}|$ such that $\eta_{C}^{ 3}\cong \mathcal O_{C}$. To this purpose, we  analyze the cup product
$$
\cup: H^0(\hol_S(-2K_{S})) \otimes H^1( \hol_S(-E)) \to H^1( \hol_S(-2K_S-E)).
$$
First we observe that there is a standard exact sequence
$$
0 \to \mathcal O_{S}(-E) \to \mathcal O_{S}(-2K_{S}-E) \to \eta_{C}^{ 3} \to 0,
$$
just because  $- C - 3K_{S} - 3L \sim -E$. Secondly, we consider the associated long exact sequence and recall that the induced map
$$
\mu_{C}: H^1(\hol_{S}(-E)) \to H^1( \hol_S(-2K_S-E))
$$
is the cup product with $s_C$. Let
$$
s_{C}^{\perp} := \ \lbrace v \in H^{1}(\hol_{S}(-E)) \ | \  v \cup s_{C} = 0 \rbrace
$$
be the $\cup$-orthogonal space of $s_{C}$. The next property is immediate.
\begin{proposition} Let $C \in |-2K_{S}|$, then  $H^{0}(\eta_{C}^{3}) =$  $s_{C}^{\perp}$. \end{proposition} \rm \par \noindent
\begin{proof}  Since $h^0(\hol_S(-2K_S-E)) = 0$, the above long exact sequence implies that  $H^{0}(\eta_{C}^{3}) = \Ker \mu_{C}$
$=$ $s_{C}^{\perp}$. \end{proof} \par \noindent
Moreover we have:
 \begin{proposition}\label{orth} Let $C \in |-2K_{S}|$ be a smooth curve. Then the following conditions are equivalent:
\begin{itemize}
\item[(1)] $\eta_{C}$ is  a non zero 3-torsion element of $\Pic^{0}(C)$;
\item[(2)] the $\cup$-orthogonal space to  $s_C$ has dimension 1;
\item[(3)] $n_{C}$ is the base locus of a pencil $Q \subset |-2K_{S}|$.
\end{itemize}
\end{proposition} 
\begin{proof} (1) $\Leftrightarrow$ (2) is clear from the previous remarks. (2) $\Rightarrow$ (3): by $\eta^{\otimes 3}(n_C) \cong 
\mathcal O_C(-2K_S)$ and $\eta^{\otimes 3}$ $\cong$ $\mathcal O_C$, it holds that $\mathcal O_C(n_C) \cong \mathcal O_C(-2K_S)$. By $h^1(\mathcal O_S)$ $=$ $0$, the restriction map $H^0(\mathcal O_S(-2K_S)) \to H^0(\mathcal O_C(-2K_S))$ is surjective. Hence $n_C$ is cut out on $C$ by a member of $\mid -2K_S \mid$. (3) $\Rightarrow$ (2) follows from reversing the argument. \end{proof}
\rm \medskip \par  \noindent

\subsection*{Schlaefli's double six}

We need to consider the six exceptional lines which are defined as follows:
$$
F_{i} := \text { strict transform of the conic through $\underline x - \lbrace x_{i} \rbrace$ under $\sigma$ }, \ (i = 1 \dots 6).
$$
The twelve lines $E_{1}, \dots, E_{6} , F_{1}, \dots, F_{6}$ form a configuration of lines on $S$, which is well known as a {\em Schlaefli' s double six}. In particular, the divisor
$$
F := F_{1} + \dots + F_{6}
$$
is the exceptional divisor of a second blow up $ \hat \sigma: S \to \mathbb P^{2}$,  defined by 
$$
|\hat L| := |5L-2E|.
$$
It is immediately checked, that $ \mathcal O_{C}(-K_{S} - \hat L) \cong \mathcal O_{C}(K_{S}+L) \cong \eta_{C}^{-1}$. Let us also point out that  $\hat \sigma | C: C \to \mathbb P^{2}$ is the morphism associated to $|\omega_{C} \otimes \eta_{C}|$. 
\begin{remark}
 For later use we observe that for $(C, \eta) \in \mathcal{R}_{4,3}$, $\Sing \Gamma_{\eta} \cup \Sing \Gamma_{\eta^{-1}}$ gives rise to a Schlaefli's double six on $S$.
\end{remark}

\noindent 
Since $h^{1}(\mathcal O_{S}(-E)) = 5$, we fix from now on the following notation:
\begin{itemize} 
\item $\mathbb P^{4} := \mathbb P (H^{1}(S, \hol_S(-E)))$,
\item $\overline v \in \mathbb P^{4}$ is the point defined by $v \in H^{1}(S, \hol_S(-E)) \setminus \{ 0\}$.
\end{itemize} 
\bigskip \par \noindent
Now we  study the $\cup$-orthogonal space $v^{\perp}$ of a vector  $v \in H^{1}(\mathcal O_{S}(-E))$. Note that $v$ corresponds to a vector of $Ext^{1}(\mathcal O_{S}(-2K_{S}), \mathcal O_{S}(-2K_{S}-E))$. Therefore a general $v$ defines an extension
\begin{equation}\label{ext}
0 \to \mathcal O_{S}(-2K_{S}-E) \to \mathcal V \to \mathcal O_{S}(-2K_{S}) \to 0,
\end{equation}
where $\mathcal V$ is a rank two  vector bundle on $S$.  It is easy to compute that
$$
\det \mathcal V \cong \mathcal O_{S}(F).
$$
Passing to the long exact sequence,  the coboundary map
$$
\partial_{v}: H^{0}(\mathcal O_{S}(-2K_{S})) \to H^{1}(\mathcal O_{S}(-E-2K_{S}))
$$
is  the cup product with $v$. Therefore it follows that  
$$
v^{\perp} = \Ker (\partial_{v}) \cong H^{0}(\mathcal V).
$$
From the same long exact sequence we see that $h^{2}(\mathcal V) = 0$. Therefore, applying Riemann-Roch to $\mathcal V$, we conclude that
$$
\dim  v^{\perp} = h^{0}(\mathcal V) \geq 2.
$$
\begin{definition} Let $ \overline v \in \mathbb P^{4}$ be the point defined by the vector $v$, then
$$
P_{\overline v} := \mathbb P (v^{\perp}) (= \mathbb{P} (H^0(\mathcal{V}))).
$$
\end{definition} \par \noindent
$P_{\overline v}$ is a linear system contained in $|-2K_{S}|$, and by the previous remarks we have
$$
\dim P_{\overline v} \geq 1.
$$
In order to describe $P_{\overline v}$, we want to understand its position with respect to the  linear subspaces $\Lambda_{i} \subset | -2K_{S}|$ which  are defined as follows:
$$
\Lambda_{i} := \lbrace C \in |-2K_{S}| \ | \ F_{i} \subset C \rbrace, \ \ i = 1 \dots 6. 
$$ 

We have the following result:

\begin{theorem} \label{secondprojection} Assume that $P_{\overline v}$ has no fixed component. Then
\begin{itemize} 
\item[(1)] $P_{\overline v}$ is a pencil,
\item[(2)] $P_{\overline v} \cap \Lambda_{i}$ is a point, for each $i = 1 \dots 6$, 
\item[(3)] the base locus of $P_{\overline v}$ is a quadratic section of $F$.
\end{itemize} 
\end{theorem}
\begin{proof} We know that $\det \mathcal V \cong \mathcal O_{S}(F)$ and $h^{0}(\mathcal V) \geq 2$. Let us consider $s \in H^{0}(\mathcal V) \setminus \{ 0 \}$ and its scheme of zeroes $Z_{s}$. Since $c_{2}(\mathcal V) = 0$, either $Z_{s} = \emptyset$ or $\dim Z_{s} = 1$. 

\noindent
\underline{\bf{Claim.}} If $P_{\overline v}$ does not have a fixed component, then there exists a $s \in H^{0}(\mathcal V) \setminus \{ 0 \}$ such that $Z_{s} = \emptyset$.

\noindent
We will prove the claim in a moment.
\noindent
We can now assume: $Z_{s} = \emptyset$. Then $s$ defines an exact sequence
$$
0 \to \mathcal O_{S} \to \mathcal V  \to \mathcal O_{S}(F) \to 0.
$$
Since $S$ is regular, the associated long exact sequence yields $ h^{0}(\mathcal V) = 2$. Hence $\mathcal P_{\overline v}$ is a pencil and (1) follows. \par \noindent  To prove (2), tensor the above exact sequence by $\mathcal O_{S}(-F_{i})$, and consider the associated long exact sequence. Since
$h^{1}(\mathcal O_{S}(-F_{i})) = 0$ and $F_{i}$ is a component of $F$, it follows 
$
h^{0}(\mathcal V(-F_{i})) = h^{0}(\mathcal O_{S}(F-F_{i})) = 1.
$ \ On the other hand, the extension defined by $v$ induces the following commutative diagram with exact lines
$$
\begin{CD}
 {H^{0}(\mathcal O_{S}(-2K_{S}-E-F_{i}))} @>>> {H^{0}(
\mathcal V(-F_{i}))} @>u>> {H^{0}(\mathcal O_{S}(-2K_{S}-F_{i}))} \\
@VVV @VVV @VVV @. \\
{H^{0}(\mathcal O_{S}(-2K_{S}-E))} @>>> {H^{0}(
\mathcal V)} @>>> {H^{0}(\mathcal O_{S}(-2K_{S}))}. \\
\end{CD}
$$
Here  $u$ is injective, because $h^{0}(\mathcal O_{S}(-2K_{S}-E-F_{i})) = 0$.  The vertical arrows are injective too. Then the image of $H^{0}(\mathcal V(-F_{i}))$ in $H^{0}(\mathcal O_{S}(-2K_{S}))$ is the unique element  $C \in P_{\overline v}$ which contains $F_{i}$. This implies (2). \par \noindent

\begin{proof}[Proof of the claim.]
Assume  that  $\dim Z_{s} = 1$, $\forall s \in H^{0}(\mathcal V) - \lbrace 0 \rbrace$. Let $s \in  H^{0}(\mathcal V) - \lbrace 0 \rbrace$ and write $Z_{s} = Z_{0,s} + D_{s}$, where $Z_{0,s}$ is 0-dimensional and $D_{s}$ is a curve. Consider the exact sequence
$$
0 \to \mathcal O_{S}(D_{s}) \to \mathcal V \to \mathcal I_{Z_{0,s}}(F-D_{s}) \to 0,
$$
where $\mathcal I_{Z_{0,s}}$ is the ideal sheaf of $Z_{0,s}$ in $S$. Applying the same argument as above just replacing $F_i$ by $D_s$, it follows that each $C \in P_{\overline v}$ contains a curve $D_{s}$. On the other hand, $D_{s} \neq C$: if not, $\mathcal O_{S}(-2K_{S})$ would be a subbundle of $\mathcal V$, which is impossible for $s \neq 0$. Thus we have shown that each $C \in P_{\overline v}$ is reducible. For
a general $C \in P_{\overline v}$,  let $d := -K_S \cdot B$ be the minimal degree of an integral curve $B \subset C$. Then one has $d \leq 3$. If $d = 1$, 
$B$ is a line and a fixed component of $P_{\overline v}$. If $d = 2$, $B$ is a conic and $dim \ \mid B \mid = 1$. It is easy to see  that the family of reducible curves $C \in \mid -2K_S \mid$ containing some conic of $\mid B \mid$ is exactly  the Segre embedding of
$\mathbf P^4 \times \mathbf P^1  :=  \mid B \mid \times \mid -2K_S - B \mid$. Note that $P_{\overline v}$ is a linear space in $\mathbf P^4 \times \mathbf P^1$ and that the projection $P_{\overline v} \to \mid B \mid$ is a surjective morphism. This implies that $P_{\overline v} = \lbrace E \rbrace \times \mid B \mid$ for some $E \in \mid -2K_S-B \mid$. Hence $P_{\overline v}$ is a pencil and $E$ is its fixed component. If $d = 3$ and $B$ is a skew cubic  then $\dim \mid B \mid = 2$. The argument is completely analogous to the previous one: we leave it to the reader. The last case is when $d = 3$ and $B$ is a plane section of $S$. Then any $C \in P_{\overline v} $ is the union of two  elements $B', B''$  of $\mid -K_S \mid$. It is well known that either one of the curves $B', B''$ is a fixed component  or there exist curves $C$ such that $B' = B''$. This implies that $\mathcal O_S(-2K_S)$ is a subbundle of $\mathcal V$, a case already excluded.  \end{proof} \par \noindent (3) The base locus $m$ of $P_{\overline v}$ has degree $12$. By (2), $P_{\overline v}$ contains some $C_{i} = D_{i} + F_{i}$,  $i = 1 \dots 6$. Let $C \in P_{\overline v}$ be an element not containing $F_{1} \dots F_{6}$. Then $m_{i} := C \cdot F_{i}$ is a subscheme of $m$. The same holds for   $C \cdot F = \sum m_{i}$, since the $F_{i}$'s are disjoint. Hence $m = C \cdot F$ for degree reasons.
 \end{proof}

\begin{definition} For any $C \in |-2K_{S}|$ we define $m_{C} := C \cdot F$. \end{definition} \noindent
The next result simply summarizes some useful equivalent conditions.
\begin{proposition} \label{cormain} Let $C \in |-2K_{S}|$ be smooth. Then the following conditions are equivalent:
\begin{itemize} 
\item[(1)] $\eta_{C}$ is a non trivial 3-torsion element of $\Pic^{0}(C)$,
\item[(2)] $C \in P_{\overline v}$ for some $\overline v \in \mathbb P^{4}$,
\item[(3)] $n_{C}$ is the base locus of a pencil $Q \subset |-2K_{S}|$,
\item[(4)] $m_{C}$ is the base locus of a pencil $P$ for some $\overline v \in \mathbb P^{4}$. 
\end{itemize} \par \noindent
Moreover, if (4) holds, then there exists $\overline v \in \mathbf P^4$ such that $P = P_{\overline v}$.
\end{proposition}
\begin{proof} Conditions (1), (2), (3) are equivalent by  proposition (\ref{orth}). To prove their equivalence with (4), we recall that $\hat L = 5L - 2E$ defines a morphism $\hat \sigma: S \to \mathbb P^{2}$, which is the contraction of the previously defined divisor $F$. It is clear that $\mathcal O_{C}(-K_{S}-\hat L) \cong \eta_{C}^{-1}$. Then, applying proposition (\ref{orth}) to $\hat \sigma$ and $\eta_{C}^{-1}$, it follows that $\eta_{C}^{-1}$ is a non trivial 3-torsion element of $\Pic^{0}(C)$ iff (4) holds.  Proposition 2.3 also implies the last statement. 
\end{proof}

\begin{proposition} \label{dimension} There exists $\overline v \in \mathbb P^{4}$, such that $P_{\overline v}$ is a pencil with no fixed component and
the general element $C \in P_{\overline v}$ is smooth.
\end{proposition}
\begin{proof} We use Schlaefli's  double six configuration $E_{1}, \dots, E_{6}, F_{1}, \dots F_{6}$.
Observe that $|-2K_{S}|$ contains the elements
$$
C_{0} := E_{1} + E_{2} + E_{3} + F_{4} + F_{5} + F_{6} \ , \ C_{\infty} := F_{1} + F_{2} + F_{3} + E_{4} + E_{5} + E_{6}.
$$
Obviously, the pencil $P$ generated by $C_{0}$ and $C_{\infty}$  has no fixed components. With a small additional effort one shows that its base locus $m$ is a smooth quadratic section of $F$. Hence a general $C \in P$ is smooth and $m = m_{C}$. Then the previous proposition (\ref{cormain}) implies that $P$ $=$ $P_{\overline v}$, for some $\overline v \in \mathbb P^{4}$.
\end{proof}
  \begin{corollary} \label{cupproduct} The cup product
$$
\cup : H^{1}(\mathcal O_{S}(-E)) \otimes H^{0}(\mathcal O_{S}(-2K_{S})) \to H^{1}(\mathcal O_{S}(-E - 2K_{S}))
$$
is surjective for any $\underline x \in \mathcal X$.
\end{corollary}
\begin{proof} Let $\cup_{v}: <v> \otimes H^{0}(\mathcal O_{S}(-2K_{S})) \to H^{1}(\mathcal O_{S}(-E-2K_{S}))$ be the cup product with
a general $v \in H^{1}(\mathcal O_{S}(-E))$. By theorem (\ref{secondprojection}) and proposition (\ref{dimension}),  $\dim \Ker \cup_{v} = v^{\perp} = 2$. Since
$\dim ( <v> \otimes H^{0}(\mathcal O_{S}(-2K_{S}))) = 10$ and  $h^{1}(\mathcal O_{S}(-E-2K_{S})) = 8$, it follows that $\cup_{v}$ is surjective. Hence $\cup$ is surjective. 
\end{proof} \noindent
We mention without proof, since it is not needed in the sequel,  a method to describe those $P_{\overline v}$'s having fixed components. 
 \begin{proposition} Let $D$ be the fixed curve of a linear system $P_{\overline v}$. Then $D$ is contained in a quadratic section of $S$. Moreover, $\mathcal V$ fits in an exact sequence
$$
0 \to \mathcal O_{S}(D) \to \mathcal V \to \mathcal I_{Z}(F-D) \to 0
$$
where $\dim Z = 0$, $\mathcal I_{Z}$ is the ideal sheaf of $Z$ and $\deg Z + D(F-D) = 0$. \end{proposition} \par \noindent
\begin{remark} For instance, let $\overline v_{i} \in \mathbb P^{4}$ be  defining the natural extension
$$
0 \to \mathcal O_{S}(-2K_{S} - E) \to  \mathcal O_{S}(-2K_{S} - E + E_{i}) \oplus \mathcal O_{S}(-2K_{S}-E_{i}) \to \mathcal O_{S}(-2K_{S}) \to  0.
$$
In this case, we have $D = E_{i}$ for the fixed curve of $P_{\overline v_{i}}$ and $Z = 0$. Moreover,
$$
P_{\overline v_{i}} = |-2K_{S}-E_{i}| = \mathbb P^{6}.
$$
\end{remark} \noindent
All this  suggests to study  the projectivized set of decomposable tensors $v \otimes s$ such that $v \cup s = 0$ and $div(s)$ is smooth. To describe it, we consider: 
\begin{itemize} 
\item $\mathbb P^{49} := \mathbb P(H^{1}(-E) \otimes H^{0}(-2K_{S}))$,
\item the Segre embedding $\mathbb P^{4} \times |-2K_{S}| \subset \mathbb P^{49}$,
\item the linear subspace $\mathbb P (\Ker \cup) \subset \mathbb P^{49}$.
\end{itemize} \noindent
Then we give the following
\begin{definition}  
$ 
\mathbb T_{\underline x} := \mathbb P(\Ker \cup) \ \cap \ \mathbb P^{4} \times |-2K_{S}|.
$
\end{definition} \noindent
It is clear from the definition that
$$
\mathbb T_{\underline x} = \lbrace (\overline v, C) \in \mathbb P^{4} \times |-2K_{S}| \ | \ h^{0}(\eta_{C}^{3}) \geq 1 \rbrace,
$$
and also that, for any $\overline v \in \mathbb P^{4}$, one has
$$
\mathbb T_{\underline x} \cap \lbrace \overline v \rbrace \times |-2K_{S}| = \lbrace \overline v \rbrace \times  P_{\overline v}.
$$
It turns out that $\mathbb T_{\underline x}$ is reducible: this is clear considering the example in remark (2.3) and theorem (\ref{fibre}). However, we are now ready to show that there exists a unique irreducible component containing all pairs $(\overline v, C) \in \mathbb T_{\underline x}$, such that $C$ is smooth. To prove this, one more definition will be convenient.
\begin{definition} \label{definition 2.2} \par \noindent
\begin{itemize} 
\item  $\mathbb M_{\underline x} :=  \text {\it Zariski closure of} \   \lbrace (\overline v, C) \in \mathbb T_{\underline x} \ | \ \dim P_{\overline v} = 1 \rbrace$;
\medskip
\item  $\mathbb M_{\underline x}^{o} := \lbrace (\overline v, C) \in \mathbb T_{\underline x} \ | \ C \ is \ smooth \rbrace$; 
\item  $p_{\underline x}: \mathbb M_{\underline x} \to \mathbb P^{4} \  and \ q_{\underline x}: \mathbb M_{\underline x} \to |-2K_{S}|$ are the projection maps.
\end{itemize} 
\end{definition} \noindent
By proposition (\ref{dimension}) $\mathbb M_{\underline x}$ and $\mathbb M_{\underline x}^{o}$ are not empty. Notice also that
$\mathbb M_{\underline x}^{o} \subset \mathbb M_{\underline x}$.

\medskip
\noindent
Indeed, let $(\overline v, C) \in \mathbb M_{\underline x}^{o}$, then $C \in P_{\overline v}$. Since $C$ is integral, $P_{\overline v}$ has no fixed component. Hence, by theorem (\ref{secondprojection}), $\dim P_{\overline v} = 1$. Clearly,  $\mathbb M_{\underline x}^{o}$ is open in $\mathbb M_{\underline x}$.

\noindent
We are now ready to conclude this section:

\begin{theorem} \label{fibre} For any $\underline x \in \mathcal X$ the projection $ p_{\underline x}: \mathbb M_{\underline x} \to \mathbb P^{4}$ is surjective. Moreover, $p_{\underline x}$ is a locally trivial $\mathbb P^{1}$-bundle  over a non empty open set
of $\mathbb P^{4}$. In particular, $\mathbb M_{\underline x}$ is irreducible and rational.
\end{theorem}
\begin{proof} The fibre of $p_{\underline x}$ at $\overline v$ is the linear space $P_{\overline v}$. By definition, $\mathbb M_{\underline x}$ is the Zariski closure of the union of the fibres of minimal dimension 1. Hence $\mathbb M_{\underline x}$ is a union of irreducible components of $\mathbb T_{\underline x}$. Each of them has dimension $\geq 5$. Indeed,  by corollary (\ref{cupproduct}), the cup product map $\cup$ is surjective. Then the codimension of $\mathbb P (\Ker \cup)$ in $\mathbb P^{49}$ is $h^{1}(\mathcal O_{S}(-2K_{S}-E)) = 8$. Hence, counting dimensions, each irreducible component of $\mathbb T_{\underline x}$ has dimension $\geq 5$. But, over a dense open set $U$ of $p_{\underline x}(\mathbb M_{\underline x})$, the fibre of $p_{\underline x}$ is $\mathbb P^{1}$. Hence
$p_{\underline x}(\mathbb M_{\underline x}) = \mathbb P^{4}$ and  $\mathbb M_{\underline x}$ is irreducible. Moreover, $\mathbb M_{\underline x}$ is a $\mathbb P^{1}$-bundle over $U$. \end{proof}

 \section {Moduli of plane models of cubic surfaces} \label{ratproof} \noindent
Starting from a pair $(C,\eta)$, we came up with  a set $\underline x \in \mathcal X$ of 6 points of $\mathbb P^{2}$ in general position.  It is now  time to globalize  our constructions over the  moduli space of  $\underline x$. This space can be viewed as the GIT-quotient $\mathcal X / \PGL(3)$ and it will be discussed in the next section. On the other hand, it is wellknown that this space is birational to another moduli space, the space of pairs defined as follows:
\begin{definition} A pair $(S, L)$ is a {\it \em plane model of a cubic surface} if
\begin{itemize}
\item $S$ is a Del Pezzo surface of degree $3$;
\item $-K_S$ is very ample;
\item $L \in \Pic(S)$ and $L^2 = 1$, $K_S\cdot L = -3$.
\end{itemize}
\end{definition} 
\noindent
\begin{definition} The moduli space of pairs $(S,L)$  will be denoted by $\mathcal P$.
\end{definition} 
\noindent
As is well known, every line bundle $L \in Pic \ (S)$ satisfying the above conditions defines a map $\sigma: S \to \mathbb P^{2}$ which is the  blow up of a set $\underline x \in \mathcal X$. It is also well known that the  assignement $(S,L) \to \underline x$ induces a  map
$$
\mathcal P \to \mathcal X / \PGL(3),
$$
which is a birational morphism. \medskip \par \noindent  For the  blow up $\sigma: S \to \mathbb P^{2}$  we will keep our usual conventions:
\begin{itemize}
\item $\underline x  = \{ x_{1} \dots x_{6} \}$ $=$ fundamental locus of $\sigma^{-1}$;
\item $E =$ the exceptional divisor of $\sigma$; 
\item $E_{i} = \sigma^{-1}(x_{i})$;  
\item $S$ is embedded in $\mathbb P^{3}$ by $|-K_{S}|$.
\end{itemize}
There are $72$ plane models of the same cubic $S$, and they come in pairs: each pair defines a Schlaefli's double six. 
\begin{remark} \rm
 As already mentioned, a smooth curve $C$ of genus $4$ together with a subgroup of order $3$ in $\Pic^0(C)$ defines a Schlaefli's double six.
\end{remark}
 \noindent Let $(S_{1}, L_{1})$ and $(S_{2},L_{2})$ be plane models of cubic surfaces, then:
\begin{definition} 
A {\em morphism of plane models of cubic surfaces} $$\psi: (S_1, L_1) \to (S_2, L_2)$$  is a morphism   $\psi: S_1 \to S_2$ such that $L_2 = \psi^* L_1$.  We will  say that
\begin{itemize}
\item [(1)] $\psi$ is an isomorphism if $\psi$ is biregular; 
 \item [(2)] $\psi$ is an automorphism if, furthermore, $(S_{1},L_{1}) = (S_{2},L_{2})$.
\end{itemize}
\end{definition} \noindent
\begin{proposition}\label{trivaut} For a general plane model $(S,L)$ of a cubic surface the only automorphism is the identity. 
\end{proposition}
\begin{proof} Let $\psi: (S, L) \to (S, L)$ be an automorphism.  The assumption  $\psi^{*}L \cong L$ implies $\psi^{*}E = E$. Then $\psi$ induces a map $\bar{\psi} \in Aut(\mathbb P^{2})$ such that  $\bar{\psi}(\underline x) = \underline x$. This implies $\bar{\psi} = id_{\mathbb{P}^2}$ and hence $\psi = id_S$.
\end{proof} \noindent
The moduli space $\mathcal P$ contains  a non empty open set
$$
\mathcal P^{o} \subset \mathcal P,
$$
which is the moduli space of pairs $(S,L)$ with trivial automorphisms group. From the general theory of
moduli of Del Pezzo surfaces,  and their explicit construction as in \cite{colombo}), it follows  that  on $\mathcal P^{o}$ there exists a universal family, representing the moduli functor. This is a pair
$$
(\mathcal S, \mathcal L)
$$
where $\mathcal S$ is  a  variety endowed with a  morphism $\pi: \mathcal S  \to \mathcal P^{o}$
and a line bundle $\mathcal L$. We call such a pair the {\it universal plane model of a cubic surface}. In particular, it has the following property:
\begin{itemize}
\item the fibre of $\pi$ at a moduli point of $(S,L)$ is  $S$,
\item the restriction of $\mathcal{L}$ to the above  fibre $S$ is $L$.
\end{itemize}  \noindent
Let $\omega$ be the relative dualizing sheaf of $\pi$. Then we consider the sheaves
\begin{itemize}
\item[(i)] $ {\mathcal H} := \pi_* ({{\mathcal \omega}^{-2}})$; 
 \item[(ii)] ${\mathcal E} := R^1\pi_*({\mathcal L}^{- 3} \otimes {\mathcal \omega}^{-1})$;
\item[(iii)] ${\mathcal F} := R^1\pi_*({\mathcal \omega}^{-3} \otimes {\mathcal L}^{-3})$.
\end{itemize} \medskip  \par   \noindent 
Let $u$ be the moduli point of $(S,L)$, for the fibre over $u$ we have: \medskip  
\begin{itemize}
\item[(i)]  $ {\mathcal H}_{u} \cong H^0(S,\hol(-2K_S))$,
\item[(ii)]  $ {\mathcal E}_{u}  \cong  H^1(S,\hol(-E))$,
\item[(iii)]  $ {\mathcal F}_{u}  \cong H^1(S,\hol(-2K_S-E))$.
\end{itemize}  \noindent
Recall  that $h^{0}(\hol_S(-2K_{S})) = 10$, $h^{1}(\hol_S(-E)) = 5$ and $h^{1}(\hol_S(-E-2K_{S})) = 8$. Therefore, by Grauert's theorem, it follows that $\mathcal H$, $\mathcal E$ and $\mathcal F$ are vector bundles:  respectively of rank 10, 5 and 8. Finally we consider the morphism 
$$
\mu: {\mathcal E} \otimes {\mathcal H} \to {\mathcal F},
$$
which fibrewise over $u$ induces the cup product
$$
\cup: H^0(S,\hol_S(-2K_S)) \otimes H^1(S,\hol_S(-E)) \to H^1(S,\hol_S(-2K_S-E)).
$$
By corollary (\ref{cupproduct}), $\cup$ is surjective. Hence $\mu$ is surjective and the kernel of 
$\mu$ is a vector bundle: we will denote it as ${\mathcal K}$. 

\noindent
We fix the following notation for the induced projective bundles: \medskip
\begin{itemize}
\item[-]  $ \mathbb K := \mathbb P({\mathcal K})$, 
\item[-]  $ \mathbb P := \mathbb P ({\mathcal E} \otimes {\mathcal H})$, 
\item[-]  $  \mathbb E := \mathbb P ({\mathcal E})$,
\item[-] $ \mathbb H := \mathbb P ({\mathcal H})$.
\end{itemize} \medskip \noindent
Since ${\mathcal F}$ has rank $8$, $\mathbb K$ has codimension $8$ in $\mathbb P$. We denote by 
$$ 
\mathbb K_{u} , \ \mathbb P_{u} , \ \mathbb E_u, \ \mathbb H_u
$$
the fibres at $u$ respectively of the projective bundles $\mathbb K$, $\mathbb P$, $\mathbb E$, $\mathbb H$. We also need to consider the set of decomposable tensors
$$
\Sigma := \mathbb E \times_{\mathcal P^{o}} \mathbb H  \subset \mathbb P.  
$$
This is the subvariety of $\mathbb P$ whose fibre at $u$ is the Segre product 
$$
\mathbb E_u \times \mathbb H_{u} \subset \mathbb P_u.
$$
In particular, $\Sigma$ is endowed with two natural projections
$$
p: \Sigma \to \mathbb E \ \text{and} \ q: \Sigma \to \mathbb H.
$$
If $u$ is the moduli point of $(S,L)$, then $\mathbb H_{u} = |-2K_{S}|$. Therefore a point of  $\mathbb H$ is a pair $(u,C)$, where $u$ is as above and $C \in |-2K_{S}|$. In particular, $\mathbb H$ contains  the following open set:
$$
\mathcal U := \lbrace (u,C) \in \mathbb H \ | \ C \ \text{is smooth} \rbrace.
$$
\begin{definition} $\mathbb M$ is the Zariski closure of 
$$
\mathbb M^{o} :=  \mathbb K \ \cap \  \mathbb E \times_{\mathcal P^{o} } \mathcal U.
$$
\end{definition} \noindent
$\mathbb M^{o}$ is a scheme over $\mathcal P^{o}$. Let $\mathbb M^{o}_{u}$ be its fibre at  $u\in \mathcal P^{o}$, then 
$$
\mathbb M^{o}_{u} =  \mathbb K_{u}  \cap \mathbb E_{u} \times U,
$$
where $U = \lbrace C \in |-2K_{S}| \ | \ C \ \text{is  smooth} \rbrace$. 

\noindent
We already met $\mathbb M^{o}_{u}$. Indeed, $u$ is the moduli point of $(S,L)$ and $|L|$ defines the  blow up $\sigma: S \to \mathbb P^{2}$ of a given set $\underline x \in \mathcal X$.  In definition (\ref{definition 2.2}) we considered
$$
\mathbb M^{o}_{\underline x} =  \mathbb P(\Ker \cup) \cap \mathbb P^{4} \times U,
$$
where $U$ is the same as above. But it is clear, that $\mathbb E_{u} = \mathbb P^{4} = \mathbb P(H^{1}(\hol_S(-E)))$ and that $\mathbb K_{u} = \mathbb P(\Ker \cup)$. Hence $\mathbb M^{o}_{u} = \mathbb M^{o}_{x}$. 

\noindent
In section 2 we have described $\mathbb M^{o}_{\underline x}$ and $\mathbb M_{\underline x}$. So we are now in position to describe $\mathbb M$ and its projection map
$$
p| \mathbb M: \mathbb M \to \mathbb E.
$$
\begin{theorem} \label {definitive} $\mathbb M$ is irreducible and dominates $\mathbb E$ via $p|\mathbb M$. This map
is a locally trivial $\mathbb P^{1}$-bundle over a non empty open set of $\mathbb E$, i. e.,
$\mathbb M$ is birational to $ \mathcal P \times \mathbb P^{4} \times  \mathbb P^{1}$.
\end{theorem} 
\begin{proof} Note that, at a general point $(u,C) \in \mathbb E$, the fibre of $p$ is 
$\mathbb M_{\underline x}$. The statement then follows from the description of $\mathbb M_{\underline x}$ given in theorem (\ref{fibre}). We omit for brevity several standard details.  \end{proof} 

\section {A birational model for $\mathcal R_{4,3}$ and the rationality of $\mathcal R_{4,\langle 3 \rangle}$} \noindent
The aim of this section is the proof  that $\mathbb M$ is birational to the moduli space $\mathcal R_{4,3}$ of \'etale triple covers of genus 4 curves. Therefore, by theorem (\ref{definitive}) we have that 
$\mathcal R_{4,3}$ is birational to $ \mathcal P \times \mathbb P^{5}$.
\begin{remark}\label{involution}
\rm  It is an obvious consequence of our construction that the involution  $i : \mathcal R_{4,3} \rightarrow \mathcal R_{4,3}$, $i(C,\eta) = (C , \eta^{-1})$ corresponds to the involution $(j,id): \mathcal P \times \mathbb P^{5} \rightarrow \mathcal P \times \mathbb P^{5}$, where $j$ is the Schlaefli involution on $\mathcal{P}$. 
\end{remark} \noindent
To begin, we  observe that a point in the open set $\mathbb M^{o}$ (defined in the previous section) is a triple 
$$(u, \overline v, C), $$ 
where
$u$ is the moduli point of a plane model of a cubic $(S,L)$ or equivalently of six unordered points $\underline{x} \in \mathcal{X}$ in general position in $\mathbb P^2$ and  $(\overline v, C) \in \mathbb M_{u} =  \mathbb M^{o}_{\underline x}$. This is equivalent to say that $C$ is a smooth element of $|-2K_{S}|$ and that
$$
v \cup s_{C} = 0,
$$
where $v \in H^{1}(\mathcal O_{S}(-E))$ defines  the point $\overline v \in \mathbb E_{u}$ and $s_{C} \in H^{0}(\mathcal O_{S}(-2K_{S}))$ is an equation of $C$.  We proved in section 2 that $\eta := \mathcal O_{C}(-L-K_{S})$ is a non
zero 3-torsion element of $\Pic^{0}(C)$. Hence $\mathbb M$ comes equipped with a natural rational map
$$
\alpha: \mathbb M \to \mathcal R_{4,3},
$$
sending the triple $(u,\overline v,C)$ to the moduli point of $(C,\eta)$. Conversely, we want now to show that the triple $(u, \overline v, C)$ is uniquely reconstructed from the pair $(C, \eta)$. This implies that
there exists a second rational map
$$
\beta: \mathcal R_{4,3} \to \mathbb M,
$$
which is inverse to $\alpha$. We proceed in several steps. \rm \par \noindent 
\em {(1)} \rm We observe that the image of $|\omega_{C} \otimes\eta^{-1} |$ is a plane sextic $\Gamma_{\eta}$, such that
$
\Sing  \Gamma_{\eta}  \equiv  \underline x  \mod \PGL(3).
$
This follows, because $\eta \cong \mathcal O_{C}(-K_{S}-L)$. Since $\omega_{C} \cong \mathcal O_{C}(-K_{S}) $, we have $|\omega_{C} \otimes \eta_{C}^{-1}| = |\mathcal O_{C}(L)|$. Hence, up to projective automorphisms,
$\Gamma_{\eta} = \sigma(C)$ and, therefore,  $\Sing \Gamma_{\eta} = \underline x \mod \PGL(3)$. \par \noindent
\em {(2) } \rm Starting from $(C,\eta)$, we first associated to it the moduli point $u$ of the pair $(S,L)$: equivalently, the moduli point of $\Sing  \Gamma_{\eta}$. As a second step, we reconstruct uniquely the curve $C$ in the linear system $|-2K_{S}|$. More precisely, we have to show that there exists a unique $D \in |-2K_{S}|$ such that $D \cong C$ and $\mathcal O_{D}(-K_{S}-L) \cong \eta$. This follows from the next  lemma.
\begin{lemma}  Assume that $C_{1}, C_{2} \in |-2K_{S}|$ are smooth biregular curves such that $\epsilon_{1} \cong \epsilon_{2}$,  where $\epsilon_{i} =\mathcal O_{C_{i}}(-K_{S}-L)$. Then $C_{1} = C_{2}$.
\end{lemma}
\begin{proof} Applying the same proof as the one of proposition (2.1),  it follows,  that $\epsilon_{i}$ is non trivial. Let $\Gamma_{i} = \sigma(C_{i})$.  Then, under our assumption, $\Gamma_{1}$ and $\Gamma_{2}$
are projectively equivalent. Since $\Sing \Gamma_{1} = \Sing \Gamma_{2} = \underline x$, there
exists $a \in \Aut(\mathbb P^{2})$ such that $a(\underline x) = \underline x$. Since the pair $(S,L)$
has no automorphism, $a$ is the identity and hence $C_{1} = C_{2}$. \end{proof} \noindent
\begin{remark} \rm Let $\Pic_{0,4}$ be the universal Picard variety of genus 4 and degree 0 and let
$s: |-2K_{S}| \to \Pic_{0,4}$ be the rational map sending $C$ to the moduli point of the pair
 $(C,\mathcal O_{C}(-K_{S}-L))$. By the previous lemma, $s$ is injective on the open set $U$ of smooth curves. \end{remark} \medskip \par  \noindent
\em {(3)} \rm  So far we reconstructed uniquely from $(C,\eta)$ the moduli point $u$ and a copy of $C$ in
$|-2K_{S}|$, such that $\eta \cong \mathcal O_{C}(-K_{S}-L)$. Finally, the point $\overline v$ is also uniquely reconstructed: consider the standard exact sequence
$$
0 \to \mathcal O_S(-E) \to \mathcal O_{S}(-E-2K_{S}) \to \eta^{3} \to 0.
$$
Passing to the long exact sequence, the image of $H^{0}(\eta^{3})$ via the coboundary map is exactly the vector space $\overline v$ generated by $v$. \medskip \par \noindent
Due to the steps (1), (2) and (3) the triple $(u, \overline v, C)$ is uniquely defined by the
isomorphism class of $(C,\eta)$. Hence there exists a rational map
$$
\beta: \mathcal R_{4,3} \to \mathbb M,
$$
sending the moduli point of $(C,\eta)$ to $(u, \overline v, C)$.  We conclude that
\begin{theorem} $\alpha: \mathbb M \to \mathcal R_{4,3}$ is birational. In particular, $\mathcal R_{4,3} \simeq \mathcal P \times \mathbb P^5$. \end{theorem}
\begin{proof} Since $\mathbb M$ and $\mathcal R_{4,3}$ are integral of the same dimension, $\alpha$ and $\beta$ are each the inverse of the other. The second part follows then from theorem (\ref{definitive}).
\end{proof}

We conclude the section by the following
\begin{theorem} $\mathcal R_{4,\langle 3 \rangle}$ is rational. 
\end{theorem}
\begin{proof} We have $\mathcal R_{4,\langle 3 \rangle} \simeq \mathcal R_{4,3} / i \simeq  \mathcal P /j \times \mathbb P^{5}$, where $j$ is the Schlaefli involution. Moreover, $\mathcal P / j$ is rational, as shown in the next section. \end{proof}

 \section {The rationality of the moduli space of double sixes in $\mathbb P^{2}$} \noindent
Recall that $\mathcal P$ is birational to the GIT-quotient $\mathcal X / PGL(3)$. In this section, the rationality of the quotient of the above moduli space of six unordered points in $\mathbb{P}^2$ by the Schlaefli involution is proved. In other words, it is shown that the moduli space of double sixes (mod $PGL(3)$) in $\mathbb{P}^2$ is rational. This is an unpublished result of Igor Dolgachev, who kindly communicated us his proof. In the remaining part of this section, we present Dolgachev's proof.
His result completes the proof of the rationality of $\mathcal R_{4,\langle 3 \rangle}$. 

\begin{theorem}[I. Dolgachev]\label{dolga}
 The moduli space of double sixes in $\mathbb{P}^2$ is rational.
\end{theorem}
\begin{proof}  Consider  the GIT-quotient $\mathcal Y := \mathcal X / PGL(3)$. Let $\Delta$ be the big
diagonal in $(\mathbb P^{2})^{6}$.  Let $P_2^6 = ((\mathbb P^{2})^{6} \setminus \Delta)/\!/\textup{SL}(3)$ be the GIT-quotient with respect to the symmetric linearization of the action of $\text{SL}(3)$. Obviously, $\mathcal Y$ is birationally isomorphic to $P_2^6/\mathfrak{S}_6$.
Recall from \cite{dolgort} that 
$$P_2^6 \cong \textup{Proj} R_2^6,$$
where $R_2^6$ is the graded algebra with $d$-homogeneous part $(R_2^6)_d = (V(d)^{\otimes 6})^{\textup{SL}(3)}$, and $V(d)=H^0(\mathbb P^2,\mathcal O_{\mathbb P^2}(d))$. 

\medskip\noindent
Let $V(1)_i$ be the $i$-th copy of the space $V(1)$ and  
let $t_1^{(i)},t_2^{(i)},t_3^{(i)}$ be a basis of $V(1)_i$. For any subset $I = \{i_1,i_2,i_3\}$ of $\{1,\ldots,6\}$, denote by $D_I$ the determinant of the matrix 
$$\begin{pmatrix} t_1^{(i_1)}&t_2^{(i_1)}&t_3^{(i_1)}\\
t_1^{(i_2)}&t_2^{(i_2)}&t_3^{(i_2)}\\
t_1^{(i_3)}&t_2^{(i_3)}&t_3^{(i_3)}\end{pmatrix}.$$
We consider $D_I$ as an element of $V(1)_{i_1}\otimes V(1)_{i_2}\otimes V(1)_{i_3}$.  By the {\em Fundamental Theorem of Invariant Theory}, the vector space $(R_2^6)_d$ is spanned by the (tensor) products of
$D_I$, such that each $k\in \{1,\ldots,6\}$ appears exactly $d$ times. An explicit computation (due to A. Coble (see \cite{dolgort})) shows, that the graded algebra $R_2^6$ is generated by $5$ elements of degree 1 given by 
\begin{multline*}
x_0 = D_{123}D_{456},\  x_1 = D_{124}D_{356}, \  x_2 = D_{125}D_{346}, \\
x_3 = D_{134}D_{256}, \ x_4 = D_{135}D_{246},
\end{multline*}
and one element of degree $2$ 
$$x_5 = D_{123}D_{145}D_{246}D_{356}-D_{124}D_{135}D_{236}D_{456}.$$
There is one relation
\begin{multline*}
x_5^2= (-x_2x_3+x_1x_4+x_0x_1+x_0x_4-x_0x_2-x_0x_3-x_0)^2 -\\
-4x_0x_1x_4(-x_0+x_1-x_2-x_3+x_4).
\end{multline*}
After change of a basis 
$$(y_0,y_1,y_2,y_3,y_4,y_5) = (x_0,x_1,x_4,-x_0-x_2,-x_0-x_3,x_5)$$
the relation becomes
\begin{equation}\label{relation}
y_5^2 = (y_0y_1+y_0y_2+y_1y_2-y_3y_4)^2-4y_0y_1y_2(y_0+y_1+y_2+y_3+y_4).
\end{equation}
Note that the polynomial on the right-hand-side defines a quartic hypersurface in $\mathbb P^4$ isomorphic to the dual of the $10$-nodal Segre cubic threefold. This hypersurface is also isomorphic to a compactification of the moduli space of principally polarized abelian surfaces with level two structure (see \cite {Igusa}). 

Now, let us see how the permutation group $\mathfrak{S}_6$ acts on the generators. Using the straightening algorithm, we find that the following representatives of the conjugacy classes of $\mathfrak{S}_6$  act  on the the space $(R_2^6)_1$, i.e., on $(x_0, \ldots , x_4)$, as follows
\begin{eqnarray*}
(12)&:& (-x_0,-x_1,-x_2,x_0-x_1+x_3,-x_0-x_2+x_4),\\
(12)(34)&:& (-x_1,-x_0,x_2,-x_0+x_1-x_3,-x_0-x_3+x_4),\\
(12)(34)(56)&:& (x_1,x_0,x_0-x_1+x_2,x_0-x_1+x_3,x_4),\\
(123)&:& (x_0,x_0-x_1+x_3,-x_0-x_2+x_4,-x_1,-x_2),\\
(1234)&:& (x_0-x_1+x_3,x_0,x_0+x_2-x_4,x_1,x_0+x_3-x_4),\\
(1234)(56)&:& (-x_0+x_1-x_3,-x_0,x_0-x_1+x_2+x_3-x_4,\\
&&-x_1,-x_4),\\
(12345)&:& (-x_0+x_1-x_3,x_0+x_2-x_4,x_0,x_0+x_3-x_4,x_1),\\
(123)(45)&:& (-x_0,-x_0-x_2+x_4,x_0-x_1+x_3,-x_2,-x_1),\\
(123456)&:& (x_0-x_1+x_3,-x_0-x_2+x_4,-x_0+x_1-x_2-x_3+x_4,\\
&&-x_0-x_3+x_4,x_4),\\
(123)(456)&:& (x_0,x_0+x_2-x_4,x_0-x_1+x_2+x_3-x_4,\\
&&x_2,x_0-x_1+x_2).
\end{eqnarray*}
 This allows to compute the characters of the representation of $\mathfrak{S}_6$ on $(R_2^6)_1$, and we find, that the above representation is isomorphic to the $5$-dimensional irreducible  $\mathbb{U}_5$, which is obtained from the $5$-dimensional {\em standard} representation of $\mathfrak{S}_6$ by composition with an outer automorphism of $\mathfrak{S}_6$.

One checks, that the polynomial
$$(y_0y_1+y_0y_2+y_1y_2-y_3y_4)^2 - 4  y_0y_1y_2(-y_0+y_1+y_2+y_3+y_4)$$
is invariant. Therefore, under the isomorphism of representations $(R_2^6)_1\cong \mathbb U_5$, this polynomial  corresponds to a linear combination of the elementary symmetric polynomials  $\sigma_2$ and $\sigma_4$ in variables $z_0,\ldots,z_4$.
Since $\mathfrak S_6$ acts on $R_2^6$ by automorphisms of the graded ring, the subspace $(R_6^2)_2$ is invariant and the relation \eqref{relation} is $\mathfrak S_6$-invariant. This shows that $\mathfrak S_6$ acts on $x_5$ by changing it to $\pm x_5$. It is immediately checked that $x_5$ is not invariant, hence it is transformed via the sign representation of $\mathfrak S_6$. There is another skew invariant of $\mathbb{U}_5$, the squareroot $D$ of the discriminant $\delta$ of a general polynomial of degree $6$, whence we get a further invariant of the $\mathfrak{S}_6$ action, namely $Dx_5$. We have the equation $(Dx_5)^2 = \delta (\lambda_1 \sigma_2^2 + \lambda_2 \sigma_4)$, which is invariant. Therefore, we see that $X_6^2$ is isomorphic to a hypersurface $F_{34} \subset \mathbb{P}(2,3,4,5,6,17)$. 

\medskip
\noindent
Now we take the quotient of $X_6^2$ by the Schlaefli involution $j$. By \cite{dolgort} the Schlaefli involution is just the association morphism $a_{2,6} : R_6^2 \rightarrow R_6^2$, which is given by
$x_i \mapsto x_i$, for $0 \leq i \leq 4$, and $D_{123}D_{145}D_{246}D_{356} \mapsto D_{124}D_{135}D_{236}D_{456}$. This implies that $x_5 \mapsto - x_5$, and the squareroot $D$ of the discriminant $\delta$ is invariant. This implies that $X_6^2 / j$ is isomorphic to $\mathbb{P}(2,3,4,5,6)$, whence rational.
\end{proof}

 \end{document}